\documentclass[12pt]{amsart}
\usepackage{amsmath,amscd,amssymb,amsfonts}
\newcommand{\ver}{Aug. 28, 2005, v.4}
\newcommand{\ssbull}{\raise.2ex\hbox{${\scriptscriptstyle\bullet}$}}
\newcommand{\scirc}{\raise.2ex\hbox{${\scriptstyle\circ}$}}
\newcommand{\mopls}{\hbox{$\bigoplus$}}
\newcommand{\msum}{\hbox{$\sum$}}
\newcommand{\mtim}{\hbox{$\times$}}

\newcommand{\bC}{{\mathbb C}}
\newcommand{\bN}{{\mathbb N}}
\newcommand{\bP}{{\mathbb P}}
\newcommand{\bQ}{{\mathbb Q}}
\newcommand{\bZ}{{\mathbb Z}}
\newcommand{\boH}{{\mathbf H}}
\newcommand{\boR}{{\mathbf R}}
\newcommand{\cB}{{\mathcal B}}
\newcommand{\cD}{{\mathcal D}}
\newcommand{\cC}{{\mathcal C}}

\newcommand{\cH}{{\mathcal H}}
\newcommand{\cI}{{\mathcal I}}
\newcommand{\cO}{{\mathcal O}}
\newcommand{\cP}{{\mathcal P}}
\newcommand{\cX}{{\mathcal X}}

\newcommand{\odel}{\overline{\delta}}
\newcommand{\tM}{\widetilde{M}}
\newcommand{\tV}{\widetilde{V}}
\newcommand{\tY}{\widetilde{Y}}
\newcommand{\CH}{\hbox{\rm CH}}
\newcommand{\Gr}{\hbox{\rm Gr}}
\newcommand{\Ext}{\hbox{\rm Ext}}
\newcommand{\Hom}{\hbox{\rm Hom}}
\newcommand{\cHom}{{\mathcal Hom}}
\newcommand{\Diff}{\hbox{\rm Diff}}
\newcommand{\DR}{\text{\rm DR}}
\newcommand{\MHS}{\text{\rm MHS}}
\newcommand{\simto}{\buildrel\sim\over\longrightarrow}

\begin{document}
\title[Direct image of logarithmic complexes]
{Direct image of logarithmic complexes\\
and infinitesimal invariants of cycles}
\author{Morihiko Saito}
\address{RIMS Kyoto University, Kyoto 606-8502 Japan}
\email{msaito@kurims.kyoto-u.ac.jp}
\date{\ver}
\begin{abstract}
We show that the direct image of the filtered logarithmic de Rham
complex is a direct sum of filtered logarithmic complexes with
coefficients in variations of Hodge structures, using a generalization
of the decomposition theorem of Beilinson, Bernstein and Deligne to
the case of filtered $D$-modules.
The advantage of using the logarithmic complexes is that we have
the strictness of the Hodge filtration by Deligne after taking
the cohomology group in the projective case.
As a corollary, we get the total infinitesimal invariant
of a (higher) cycle in a direct sum of the cohomology of filtered
logarithmic complexes with coefficients, and this is essentially
equivalent to the cohomology class of the cycle.
\end{abstract}
\maketitle

\bigskip
\centerline{Introduction}

\bigskip\noindent
Let
$ X $,
$ S $ be complex manifolds or
smooth algebraic varieties over a field of characteristic zero.
Let
$ f : X \to S $ be a projective morphism, and
$ D $ be a divisor on
$ S $ such that
$ f $ is smooth over
$ S \setminus D $.
We have a filtered locally free
$ \cO $-module
$ (V^{i},F) $ on
$ S \setminus D $ underlying a variation of Hodge structure
whose fiber
$ V_{s}^{i} $ at
$ s \in S \setminus D $ is the cohomology of the fiber
$ H^{i}(X_{s},\bC) $.
If
$ D $ is a divisor with normal crossings on
$ S $,
let
$ \tV^{i} $ denote the Deligne extension [7] of
$ V^{i} $ such that the the eigenvalues of the residue of
the connection are contained in
$ [0,1) $.
The Hodge filtration
$ F $ is naturally extended to
$ \tV^{i} $ by [25].
We have the logarithmic de Rham complex
$$
\DR_{\log}(\tV^{i}) =
\Omega_{S}^{\ssbull}(\log D)\otimes_{\cO} \tV^{i},
$$
which has the Hodge filtration
$ F^{p} $ defined by
$ \Omega_{S}^{j}(\log D)\otimes_{\cO}F^{p-j}\tV^{i} $.
In general,
$ V^{i} $ can be extended to a regular holonomic
$ \cD_{S} $-module
$ M^{i} $ on which a local defining equation of
$ D $ acts bijectively.
By [23],
$ M^{i} $ and hence the de Rham complex
$ \DR(M^{i}) $ have the Hodge filtration
$ F $.
If
$ Y := f^{*}D $ is a divisor with normal crossings on
$ X, $ then
$ \Omega_{X}^{\ssbull}(\log Y) $ has the Hodge filtration
$ F $ defined by the truncation
$ \sigma $ (see [8]) as usual, i.e.
$ F^{p}\Omega_{X}^{\ssbull}(\log Y) =
\Omega_{X}^{\ssbull\ge p}(\log Y) $.

\medskip\noindent
{\bf Theorem~1.}
{\it Assume
$ Y = f^{*}D $ is a divisor with normal crossings.
There is an increasing split filtration
$ L $ on the filtered complex
$ \boR f_{*}(\Omega_{X}^{\ssbull}(\log Y),F) $ such that we have
noncanonical and canonical isomorphisms in the filtered
derived category:
$$
\aligned
\boR f_{*}(\Omega_{X}^{\ssbull}(\log Y),F)
&\simeq\mopls_{i\in \bZ}(\DR(M^{i}),F)[-i],
\\
\Gr_{i}^{L}\boR f_{*}(\Omega_{X}^{\ssbull}(\log Y),F)
&=(\DR(M^{i}),F)[-i].
\endaligned
$$
If
$ D $ is a divisor with normal crossings, we have also
}
$$
\aligned
\boR f_{*}(\Omega_{X}^{\ssbull}(\log Y),F)
&\simeq\mopls_{i\in \bZ}(\DR_{\log}(\tV^{i}),F)[-i],
\\
\Gr_{i}^{L}\boR f_{*}(\Omega_{X}^{\ssbull}(\log Y),F)
&=(\DR_{\log}(\tV^{i}),F)[-i].
\endaligned
$$

\medskip
This follows from the decomposition theorem (see [2])
extended to the case of the direct image of
$ (\cO_{X},F) $ as a filtered
$ \cD $-module, see [22].
Note that Hodge modules do not appear in the last statement if
$ D $ is a divisor with normal crossings.
The assertion becomes more complicated in the non logarithmic case,
see Remark~(i) in (2.5).
A splitting of the filtration
$ L $ is given by choosing the first noncanonical isomorphism in
the filtered decomposition theorem, see (1.4.2).
A canonical choice of the splitting is given
by choosing an relatively ample class, see [9].

Let
$ \CH^{p}(X\setminus Y,n) $ be Bloch's higher Chow group, see [3].
In the analytic case, we assume for simplicity that
$ f : (X,Y)\to (S,D) $ is the base change of a projective morphism
of smooth complex algebraic varieties
$ f' : (X',Y')\to (S',D') $ by an open embedding of complex manifolds
$ S \to S'_{\mathrm an} $, and an element of
$ \CH^{p}(X\setminus Y,n) $ is the restriction of an element of
$ \CH^{p}(X'\setminus Y',n) $ to
$ X\setminus Y $.
If
$ n = 0 $, we may assume that it is the restriction of
an analytic cycle of codimension
$ p $ on
$ X $.
From Theorem~1, we can deduce

\medskip\noindent
{\bf Corollary 1.} {\it
With the above notation and assumption, let
$ \xi \in \CH^{p}(X\setminus Y,n) $.
Then, choosing a splitting of the filtration
$ L $ in Theorem~{\rm 1 (}or more precisely, choosing the first
noncanonical isomorphism in the filtered decomposition theorem
{\rm (1.4.2))}, we have the total infinitesimal invariant
$$
\aligned
&\delta_{S,D}(\xi) = (\delta_{S,D}^{i}(\xi)) \in
\mopls_{i\ge 0}\boH^{i}(S,F^{p}\DR(M^{2p-n-i})),
\\
(\text{resp.}\,\,\,
&\odel_{S,D}(\xi) = (\odel_{S,D}^{i}(\xi))\in
\mopls_{i\ge 0}\boH^{i}(S,\Gr_{F}^{p}\DR(M^{2p-n-i})),)
\endaligned
$$
where
$ \delta_{S,D}^{i}(\xi) $ {\rm (}resp.
$ \odel_{S,D}^{i}(\xi)) $ is independent of the choice of a splitting
if the
$ \delta_{S,D}^{j}(\xi) $ {\rm (}resp.
$ \odel_{S,D}^{j}(\xi)) $ vanish for
$ j < i $.
In the case
$ D $ is a divisor with normal crossings, the assertion holds with
$ \DR(M^{2p-n-i}) $ replaced by
$ \DR_{\log}(\tV^{2p-n-i}) $.
}

\medskip
This shows that the infinitesimal invariants in
[14], [13], [27], [5], [1], [24] can be
defined naturally in the cohomology of filtered logarithmic
complexes with coefficients in variations of Hodge structures if
$ D $ is a divisor with normal crossings,
see (2.4) for the compatibility with [1].
Note that if
$ S $ is Stein or affine, then
$ \boH^{i}(S,F^{p}\DR_{\log}(\tV^{q})) $ is the
$ i $-th cohomology group of the complex whose
$ j $-th component is
$ \Gamma (S,\Omega_{S}^{j}(\log D)\otimes_{\cO}F^{p-j}\tV^{q}) $.
If
$ D $ is empty, then an inductive definition of
$ \delta_{S,D}^{i}(\xi) $,
$ \odel_{S,D}^{i}(\xi) $ was given by Shuji Saito [24]
using the filtered Leray spectral sequence together with the
$ E_{2} $-degeneration argument in [6].
He also showed that the infinitesimal invariants depend only on the
cohomology class of the cycle.
If
$ S $ is projective, then it follows from [8] that the total
infinitesimal invariant
$ (\delta_{S,D}^{i}(\xi)) $ is equivalent to the cycle class of
$ \xi $ in
$ H_{\DR}^{2p-n}(X\setminus Y) $ by the strictness of the Hodge
filtration, and the filtration
$ L $ comes from the Leray filtration on the cohomology of
$ X \setminus Y $, see Remark~(iii) in (2.5).

Corollary~1 is useful to study the behavior of the
infinitesimal invariants near the boundary of the variety.
If
$ D $ is empty, let
$ \delta_{S}^{i}(\xi) $ denote
$ \delta_{S,D}^{i}(\xi) $.
We can define
$ \delta_{\DR,S}^{i}(\xi) $ as in [19] by omitting
$ F^{p} $ before
$ \DR $ in Corollary~1 where
$ D = \emptyset $.

\medskip\noindent
{\bf Corollary 2.}
{\it Assume
$ S $ is projective.
Let
$ U = S \setminus D $.
Then for each
$ i \ge 0 $,
$ \delta_{S,D}^{i}(\xi) $,
$ \odel_{S,D}^{i}(\xi) $,
$ \delta_{U}^{i}(\xi) $ and
$ \delta_{\DR,U}^{i}(\xi) $ are equivalent to each other, i.e.
one of them vanishes if and only if the others do.
}

\medskip
Indeed,
$ (\delta_{\DR,U}^{i}(\xi)) $ is determined by
$ (\delta_{U}^{i}(\xi)) $, and
$ (\delta_{U}^{i}(\xi)) $ by
$ (\delta_{S,D}^{i}(\xi)) $.
Moreover,
$ (\delta_{S,D}^{i}(\xi)) $ is equivalent to
$ (\delta_{\DR,U}^{i}(\xi)) $ by the strictness of the Hodge
filtration [8] applied to
$ (X, Y) $ together with Theorem~1, see (2.3).
For the relation with
$ \odel_{S,D}^{i}(\xi) $, see (2.1).
Note that the equivalence between
$ \delta_{U}^{i}(\xi) $ and
$ \delta_{\DR,U}^{i}(\xi) $ in the case of algebraic cycles (i.e.
$ n = 0) $ was first found by J.D.~Lewis and Shuji Saito
in [19] (assuming a conjecture of Brylinski and Zucker and the Hodge
conjecture and using an $ L^{2} $-argument).
The above arguments seem to be closely related with their question,
see also Remark~(i) in (2.5) below.

As another corollary of Theorem~1 we have

\medskip\noindent
{\bf Corollary~3.} {\it
Assume
$ f $ induces an isomorphism over
$ S \setminus D $, and
$ Y = f^{*}D $ is a divisor with normal crossings on
$ X $.
Then
}
$$
R^{i}f_{*}\Omega_{X}^{p}(\log Y) = 0\quad \text{if}\,\,i + p > \dim X.
$$

\medskip
This follows immediately from Theorem~1 since
$ M^{i} = 0 $ for
$ i \ne 0 $.
Corollary~3 is an analogue of the vanishing theorem of
Kodaira-Nakano.
However, this does not hold for a non logarithmic complex
(e.g. if
$ f $ is a blow-up with a point center).
This corollary was inspired by a question of A.~Dimca.

I would like to thank Dimca, Lewis and Shuji Saito for good questions
and useful suggestions.

In Section~1, we prove Theorem~1 after reviewing some basic facts on
filtered differential complexes.
In Section~2 we explain the application of Theorem~1 to the
infinitesimal invariants of (higher) cycles.
In Section~3 we give some examples using Lefschetz pencils.

\bigskip\bigskip
\centerline{{\bf 1. Direct image of logarithmic complexes}}

\bigskip\noindent
{\bf 1.1.~Filtered differential complexes.}
Let
$ X $ be a complex manifold or a smooth algebraic variety over a field
of characteristic zero.
Let
$ D^{b}F(\cD_{X}) $ (resp.
$ D^{b}F(\cD_{X})^{r}) $ be the bounded derived category
of filtered left (resp. right)
$ \cD_{X} $-modules.
Let
$ D^{b}F(\cO_{X},\Diff) $ be the bounded derived category of
filtered differential complexes
$ (L,F) $ where
$ F $ is exhaustive and locally bounded below (i.e.
$ F_{p} = 0 $ for
$ p \ll 0 $ locally on
$ X $), see [22], 2.2.
We have an equivalence of categories
$$
\DR^{-1} : D^{b}F(\cO_{X},\Diff) \to D^{b}F(\cD_{X})^{r},
\leqno(1.1.1)
$$
whose quasi-inverse is given by the de Rham functor
$ \DR^{r} $ for right
$ \cD $-modules, see (1.2) below.
Recall that, for a filtered
$ \cO_{X} $-module
$ (L,F) $,
the associated filtered right
$ \cD $-module
$ \DR^{-1}(L,F) $ is defined by
$$
\DR^{-1}(L,F) = (L,F)\otimes_{\cO}(\cD,F),
\leqno(1.1.2)
$$
and the morphisms
$ (L,F) \to (L',F) $ in
$ MF(\cO_{X},\Diff) $ correspond bijectively to the morphisms of
filtered
$ \cD $-modules
$ \DR^{-1}(L,F) \to \DR^{-1}(L',F) $.
More precisely, the condition on
$ (L,F) \to (L',F) $ is that the composition
$$
F_{p}L \to L \to L' \to L'/F_{q}L'
$$
is a differential operator of order
$ \le p - q - 1 $.
The proof of (1.1.1) can be reduced to the canonical filtered
quasi-isomorphism for a filtered right
$ \cD $-module
$ (M,F) $
$$
\DR^{-1}\scirc\DR^{r}(M,F) \to (M,F),
$$
which follows from a calculation of a Koszul complex.

Note that the direct image
$ f_{*} $ of filtered differential complexes is defined by the
sheaf-theoretic direct image
$ \boR f_{*} $,
and this direct image is compatible with the direct image
$ f_{*} $ of filtered
$ \cD $-modules via (1.1.1), see [22], 2.3.
So we get
$$
\boR f_{*} = \DR^{r}\scirc f_{*}\scirc \DR^{-1} : D^{b}
F(\cO_{X},\Diff) \to D^{b}F(\cO_{S},\Diff),
\leqno(1.1.3)
$$
where we use
$ \DR^{r} $ for right
$ \cD $-modules (otherwise there is a shift of complex).

\medskip\noindent
{\bf 1.2.~De Rham complex.}
The de Rham complex
$ \DR^{r}(M,F) $ of a filtered right
$ \cD $-module
$ (M,F) $ is defined by
$$
(\DR^{r}(M,F))^{i} = \hbox{$\bigwedge$}^{-i}\Theta_{X}\otimes_{\cO}
(M,F[-i])\quad \text{for}\,\,i \le 0.
\leqno(1.2.1)
$$
Here
$ (F[-i])_{p} = F_{p+i} $ in a compatible way with
$ (F[-i])^{p} = F^{p-i} $ and
$ F_{p} = F^{-p} $.
Recall that the filtered right
$ \cD $-module associated with a filtered left
$ \cD $-module
$ (M,F) $ is defined by
$$
(M,F)^{r} := (\Omega_{X}^{\dim X},F)\otimes_{\cO}(M,
F),
\leqno(1.2.2)
$$
where
$ \Gr_{p}^{F}\Omega_{X}^{\dim X} = 0 $ for
$ p \ne -\dim X $.
This induces an equivalence of categories between the left and
right
$ \cD $-modules. The usual de Rham complex
$ \DR(M,F) $ for a left
$ \cD $-module is defined by
$$
(\DR(M,F))^{i} = \Omega_{X}^{i}\otimes_{\cO}(M,
F[-i])\quad \text{for}\,\,i \ge 0,
\leqno(1.2.3)
$$
and this is compatible with (1.2.1) via (1.2.2) up to a shift
of complex, i.e.
$$
\DR(M,F) = \DR^{r}(M,F)^{r}[-\dim X].
\leqno(1.2.4)
$$

\medskip\noindent
{\bf 1.3.~Logarithmic complex.}
Let
$ X $ be as in (1.1), and
$ Y $ be a divisor with normal crossings on
$ X $.
Let
$ (V,F) $ be a filtered locally free
$ \cO $-module underlying a polarizable variation of Hodge
structure on
$ X \setminus Y $.
Let
$ (\tV,F) $ be the Deligne extension of
$ (V,F) $ to
$ X $ such that the eigenvalues of the residue of the
connection are contained in
$ [0,1) $.
Then we have the filtered logarithmic de Rham complex
$ \DR_{\log}(\tV,F) $ such that
$ F^{p} $ of its
$ i $-th component is
$$
\Omega_{X}^{i}(\log Y)\otimes F^{p-i}\tV.
$$
If
$ (M,F) = (\cO_{X},F) $ with
$ \Gr_{p}^{F}\cO_{X} = 0 $ for
$ p \ne 0 $,
then
$$
\DR_{\log}(\cO_{X},F) = (\Omega_{X}^{\ssbull}
(\log X),F).
$$

Let
$ \tV(*Y) $ be the localization of
$ \tV $ by a local defining equation of
$ Y $.
This is a regular holonomic left
$ \cD_{X} $-module underlying a mixed Hodge module, and has
the Hodge filtration
$ F $ which is generated by the Hodge filtration
$ F $ on
$ \tV $,
i.e.
$$
F_{p}\tV(*Y) = \msum_{\nu} \partial^{\nu}F^{-p+|\nu |}\tV,
$$
where
$ F_{p} = F^{-p} $ and
$ \partial^{\nu} = \prod_{i} \partial_{i}^{\nu_{i}} $ with
$ \partial_{i} = \partial /\partial x_{i} $.
Here
$ (x_{1}, \dots , x_{n}) $ is a local coordinate system such that
$ Y $ is contained in
$ \{x_{1}\cdots x_{n} = 0\} $.
By [23], 3.11, we have a filtered quasi-isomorphism
$$
\DR_{\log}(\tV,F) \simto \DR(\tV(*Y),F).
\leqno(1.3.1)
$$
This generalizes the filtered quasi-isomorphism in [7]
$$
(\Omega_{X}^{\ssbull}(\log Y),F) \simto \DR(\cO_{X}(*Y),F).
\leqno(1.3.2)
$$

Note that the direct image of the filtered
$ \cD_{X} $-module
$ (\tV(*Y),F) $ by
$ X \to pt $ in the case
$ X $ projective (or proper algebraic) is given by the cohomology
group of the de Rham complex
$ \DR(\tV(*Y),F) $ (up to a shift of complex) by definition,
and the Hodge filtration
$ F $ on the direct image is strict by the theory of Hodge modules.
So we get
$$
\aligned
F^{p}\boH^{i}(X\setminus Y,\DR(V))
&:=
F^{p}\boH^{i}(X,\DR(\tV(*Y)))
\\
=\boH^{i}(X,F^{p}\DR(\tV(*Y)))
&=
\boH^{i}(X,F^{p}\DR_{\log}(\tV)).
\endaligned
\leqno(1.3.3)
$$

\medskip\noindent
{\bf 1.4.~Decomposition theorem.}
Let
$ f : X \to S $ be a projective morphism of complex manifolds or
smooth algebraic varieties over a field of characteristic zero.
Then the decomposition theorem of Beilinson, Bernstein and Deligne
[2] is extended to the case of Hodge modules ([22], [23]),
and we have noncanonical and canonical isomorphisms
$$
\aligned
f_{*}(\cO_{X},F)
\simeq \mopls_{j\in \bZ} \cH^{j}f_{*}(\cO_{X},F)[-j]
&\quad \text{in}\,\,D^{b}F(\cD_{S}),
\\
\cH^{j}f_{*}(\cO_{X},F)
= \mopls_{Z\subset S} (M_{Z}^{j},F)
&\quad \text{in}\,\,MF(\cD_{S}),
\endaligned
\leqno(1.4.1)
$$
where
$ Z $ are irreducible closed analytic or algebraic subsets of
$ S $, and
$ (M_{Z}^{j},F) $ are filtered
$ \cD_{S} $-modules underlying a pure Hodge module of weight
$ j + \dim X $ and with strict support
$ Z $, i.e.
$ M_{Z}^{j} $ has no nontrivial sub nor quotient module
whose support is strictly smaller than
$ Z $.
(Here
$ MF(\cD_{S}) $ denotes the category of filtered left
$ \cD_{S} $-modules.)
Indeed, the second canonical isomorphism follows from the
strict support decomposition which is part of the definition of
pure Hodge modules, see [22], 5.1.6.
The first noncanonical isomorphism follows from the strictness of
the Hodge filtration and the relative hard Lefschetz theorem for
the direct image (see [22], 5.3.1) using the
$ E_{2} $-degeneration argument in [6] together with the
equivalence of categories
$ D^{b}F(\cD_{S}) \simeq D^{b}G(\cB_{S}) $.
Here
$ \cB_{S} = \mopls_{i\in\bN}F_{i}\cD_{S} $ and
$ D^{b}G(\cB_{S}) $ is the derived category of bounded complexes
of graded left
$ \cB_{S} $-modules
$ M_{\ssbull}^{\ssbull} $ such that
$ M_{i}^{j} = 0 $ for
$ i \ll 0 $ or
$ |j| \gg 0 $, see [22], 2.1.12.
We need a derived category associated to some
abelian category in order to apply the argument in [6]
(see also [9]).
In the algebraic case, we can also apply [6] to the derived
category of mixed Hodge modules on
$ S $ and it is also possible to use [23], 4.5.4 to show the
first noncanonical isomorphism.

If
$ f $ is smooth over the complement of a divisor
$ D \subset S $ and
$ Y := f^{*}D $ is a divisor with normal crossings,
then the filtered direct image
$ f_{*}(\cO_{X}(*Y),F) $ is strict (see [23], 2.15),
and we have noncanonical and canonical isomorphisms
$$
\aligned
f_{*}(\cO_{X}(*Y),F)
\simeq \mopls_{j\in \bZ} \cH^{j}f_{*}(\cO_{X}(*Y),F)[-j]
&\quad \text{in}\,\,D^{b}F(\cD_{S}),
\\
\cH^{j}f_{*}(\cO_{X}(*Y),F)
= (M_{S}^{j}(*D),F)
&\quad \text{in}\,\,MF(\cD_{S}).
\endaligned
\leqno(1.4.2)
$$
Here
$ (M_{S}^{j}(*D),F) $ is the `localization' of
$ (M_{S}^{j},F) $ along
$ D $ which is the direct image of
$ (M_{S}^{j},F)|_{U} $ by the open embedding
$ U := S \setminus D \to S $ in the category of filtered
$ \cD $-modules underlying mixed Hodge modules.
(By the Riemann-Hilbert correspondence, this gives the direct image
in the category of complexes with constructible cohomology because
$ D $ is a divisor.)
The Hodge filtration
$ F $ on the direct image is determined by using the
$ V $-filtration of Kashiwara and Malgrange, and
$ (M_{S}^{j}(*D),F) $ is the unique extension of
$ (M_{S}^{j},F)|_{U} $ which underlies a mixed
Hodge module on
$ S $ and whose underlying
$ \cD_{S} $-module is the direct image in the category of
regular holonomic
$ \cD_{S} $-modules, see [23], 2.11.
So the second canonical isomorphism follows because the left-hand
side satisfies these conditions.
(Note that
$ (M_{Z}^{j},F) $ for
$ Z \ne S $ vanishes by the localization, because
$ Z \subset D $ if
$ (M_{Z}^{j},F) \ne 0 $ .)
The first noncanonical isomorphism follows from the strictness
of the Hodge filtration and the relative hard Lefschetz theorem by
the same argument as above.

\medskip\noindent
{\bf 1.5.~Proof of Theorem~1.}
Let
$ r = \dim X - \dim S $.
By (1.1.3), (1.3.2) and (1.4.2), we have isomorphisms
$$
\aligned
\boR f_{*}(\Omega_{X}^{\ssbull}(\log Y),F)
&= \DR^{r}\scirc f_{*}\scirc \DR^{-1}(\Omega_{X}^{\ssbull}(\log Y),F)
\\
&= \DR^{r}\scirc f_{*}(\cO_{X}(*Y),F)[-\dim X]
\\
&\simeq \mopls_{i\in \bZ} \DR(M_{S}^{i}(*D),F)[-r-i],
\endaligned
\leqno(1.5.1)
$$
where the shift of complex by
$ r $ follows from the difference of the de Rham complex for
left and right
$ \cD $-modules.
Furthermore, letting
$ L $ be the filtration induced by
$ \tau $ on the complex of filtered
$ \cD_{S} $-modules
$ f_{*}(\cO_{X}(*Y),F)[-r] $,
we have a canonical isomorphism
$$
\Gr_{i}^{L}f_{*}(\cO_{X}(*Y),F)[-r] = (M_{S}^{i-r}(*D),F)[-i],
\leqno(1.5.2)
$$
and the first assertion follows by setting
$ M^{i} = M_{S}^{i-r}(*D) $.
The second assertion follows from the first by (1.3.1).
This competes the proof of Theorem~1.

\bigskip\bigskip
\centerline{{\bf 2. Infinitesimal invariants of cycles}}

\bigskip\noindent
{\bf 2.1.~Cycle classes.}
Let
$ X $ be a complex manifold, and
$ \cC^{\ssbull,\ssbull} $ denote the double complex of
vector spaces of currents on
$ X $.
The associated single complex is denoted by
$ \cC^{\ssbull} $.
Let
$ F $ be the Hodge filtration by the first index of
$ \cC^{\ssbull,\ssbull} $ (using the truncation
$ \sigma $ in [8]).
Let
$ \xi $ be an analytic cycle of codimension
$ p $ on
$ X $.
Then it is well known that
$ \xi $ defines a closed current in
$ F^{p}\cC^{2p} $ by integrating the restrictions of
$ C^{\infty} $ forms with compact supports on
$ X $ to the smooth part of the support of
$ \xi $ (and using a triangulization or a resolution of singularities
of the cycle).
So we have a cycle class of
$ \xi $ in
$ H^{2p}(X,F^{p}\Omega_{X}^{\ssbull}) $.

Assume
$ X $ is a smooth algebraic variety over a field
$ k $ of characteristic
zero.
Then the last assertion still holds (where
$ \Omega_{X}^{\ssbull} $ means
$ \Omega_{X/k}^{\ssbull} $), see [11].
Moreover, for the higher Chow groups, we have the cycle map
(see [4], [10], [12], [15], [16])
$$
cl : \CH^{p}(X,n) \to F^{p}H_{\DR}^{2p-n}(X),
$$
where the Hodge filtration
$ F $ is defined by using a smooth compactification of
$ X $ whose complement is a divisor with normal crossings,
see [8].
This cycle map is essentially equivalent to the cycle map to
$ \Gr_{F}^{p}H_{\DR}^{2p-n}(X) $ because we can reduce to the
case
$ k = \bC $ where we have the cycle map
$$
cl : \CH^{p}(X,n) \to \Hom_{\MHS}(\bQ,H^{2p-n}(X,\bQ)(p)),
$$
and morphisms of mixed Hodge structures are strictly
compatible with the Hodge filtration
$ F $.

\medskip\noindent
{\bf 2.2.~Proof of Corollary~1.}
By (2.1)
$ \xi $ has the cycle class in
$$
H^{2p-n}(X,F^{p}\Omega_{X}^{\ssbull}(\log Y)).
$$
By theorem 1, this gives the total infinitesimal invariant
$$
\delta_{S,D}(\xi) = (\delta_{S,D}^{2p-n-i}(\xi)) \in
\mopls_{i\in \bZ}\boH^{2p-n-i}(S,F^{p}\DR(M^{i})),
$$
and similarly for
$ \odel_{S,D}(\xi) $.
So the assertion follows.

\medskip\noindent
{\bf 2.3.~Proof of Corollary 2.}
Choosing the first noncanonical isomorphism in the filtered
decomposition theorem (1.4.2), we get canonical morphisms compatible
with the direct sum decompositions
$$
\aligned
\mopls_{i\ge 0}\boH^{i}(S,F^{p}\DR(M^{q-i}))
&\to
\mopls_{i\ge 0}\boH^{i}(S\setminus D,F^{p}\DR(M^{q-i}))
\\
&\to
\mopls_{i\ge 0}\boH^{i}(S\setminus D,\DR(M^{q-i})),
\endaligned
$$
and these are identified with the canonical morphisms
$$
\aligned
\boH^{q}(X,F^{p}\Omega_{X}^{\ssbull}(\log Y))
&\to
\boH^{q}(X\setminus Y,F^{p}\Omega_{X\setminus Y}^{\ssbull})
\\
&\to
\boH^{q}(X\setminus Y,\Omega_{X\setminus Y}^{\ssbull}).
\endaligned
$$
By Deligne [8], the composition of the last two morphisms is injective
because of the strictness of the Hodge filtration,
see also (1.3).
So we get the equivalence of
$ \delta_{S,D}^{i}(\xi) $,
$ \delta_{U}^{i}(\xi) $,
$ \delta_{\DR,U}^{i}(\xi) $.
The equivalence with
$ \odel_{S,D}^{i}(\xi) $ follows from (2.1).

\medskip\noindent
{\bf 2.4.~Compatibility with the definition in [1].}
When
$ D $ is empty, the infinitesimal invariants are defined in [1] by
using the extension groups of filtered
$ \cD $-modules together with the forgetful functor from the
category of mixed Hodge modules to that of filtered
$ \cD $-modules.
Its compatibility with the definition in this paper follows from
the equivalence of categories (1.1.1) and the compatibility of
the direct image functors (1.1.3).

Note that for
$ (L,F) \in D^{b}F(\cO_{X},\Diff) $ in the notation of (1.1),
we have a canonical isomorphism
$$
\Ext^{i}((\Omega_{X}^{\ssbull},F),(L,F)) = \boH^{i}(X,F_{0}L),
\leqno(2.4.1)
$$
where the extension group is taken in
$ D^{b}F(\cO_{X},\Diff) $.
Indeed, the left-hand side is canonically isomorphic to
$$
\aligned
&\Ext^{i}(\DR^{-1}(\Omega_{X}^{\ssbull},F),\DR^{-1}(L,F))
\\
&\,\,=\boH^{i}(X,F_{0}\cHom_{\cD}(\DR(\cD_{X},F),\DR^{-1}(L,F))),
\\
&\,\,=\boH^{i}(X,F_{0}\DR^{r}\DR^{-1}(L)),
\endaligned
$$
and the last group is isomorphic to the right-hand side
of (2.4.1) which is independent of a representative of
$ (L,F) $.
If
$ X $ is projective, then this assertion follows also from the adjoint
relation for filtered
$ \cD $-modules.

If
$ X $ is smooth projective and
$ Y $ is a divisor with normal crossings, then the cycle class can be
defined in
$$
\aligned
\Ext^{2p}((\Omega_{X}^{\ssbull},F),\Omega_{X}^{\ssbull}(\log Y),F[p]))
&= \boH^{2p}(X,F^{p}\Omega_{X}^{\ssbull}(\log Y))
\\
&= F^{p}\boH^{2p}(X,\Omega_{X}^{\ssbull}(\log Y)).
\endaligned
$$

\medskip\noindent
{\bf 2.5.~Remarks.}
(i) If we use (1.4.1) instead of (1.4.2) we get an analogue of
Theorem~1 for non logarithmic complexes.
However, the assertion becomes more complicated, and we get
noncanonical and canonical isomorphisms
$$
\aligned
\boR f_{*}(\Omega_{X}^{\ssbull},F)
&\simeq\mopls_{i\in\bZ,Z\subset S} (\DR(M_{Z}^{i-r}),F)[-i].
\\
\Gr_{i}^{L}\boR f_{*}(\Omega_{X}^{\ssbull},F)
&=\mopls_{Z\subset S} (\DR(M_{Z}^{i-r}),F)[-i].
\endaligned
\leqno(2.5.1)
$$
This implies an analogue of Corollary~1.
If
$ D $ is a divisor with normal crossings,
we have a filtered quasi-isomorphism for
$ Z = S $
$$
(\DR_{\log}(\tM_{S}^{i-r}),F) \simto
(\DR(M_{S}^{i-r}),F),
\leqno(2.5.2)
$$
where
$ \DR_{\log}(\tM_{S}^{i-r}) $ is the intersection of
$ \DR(M_{S}^{i-r}) $ with
$ \DR_{\log}(\tV_{S}^{i}) $.
This seems to be related with a question of Lewis and Shuji Saito,
see also [19].

\medskip
(ii) If
$ \dim S = 1 $, we can inductively define the infinitesimal invariants
in Corollary~1 by an argument similar to [24] using [26].

\medskip
(iii)
Assume
$ S $ is projective and
$ D $ is a divisor with normal crossings.
Then the Leray filtration for
$ X \to S \to pt $ is given by the truncation
$ \tau $ on the complex of filtered
$ \cD_{S} $-modules
$ f_{*}(\cO_{X}(*Y),F) $, and gives the
Leray filtration on the cohomology of
$ X \setminus Y $ (induced by the truncation
$ \tau $ as in [8]).
Indeed, the graded pieces
$ \cH^{j}f_{*}(\cO_{X}(*Y),F) $ of the filtration
$ \tau $ on
$ S $ coincide with
$ (\tV^{j+r}(*D), F) $, and give the open direct images by
$ U \to S $ of the graded pieces
$ (V^{j+r}, F) $ of the filtration
$ \tau $ on
$ U $ as filtered
$ \cD $-modules underlying mixed Hodge modules.
Note that the morphism
$ U \to S $ is open affine so that the direct image preserves
regular holonomic
$ \cD $-modules.

\newpage
\centerline{{\bf 3. Examples}}

\bigskip\noindent
{\bf 3.1.~Lefschetz pencils.}
Let
$ Y $ be a smooth irreducible projective variety of dimension
$ n $ embedded in a projective space
$ \cP $ over
$ \bC $.
We assume that
$ Y \ne \cP $ and
$ Y $ is not contained in a hyperplane of
$ \cP $ so that the hyperplane sections of
$ Y $ are parametrized by the dual projective spaces
$ \cP^{\vee} $.
Let
$ D \subset \cP^{\vee} $ denote the discriminant.
This is the image of a projective bundle over
$ Y $ (consisting of hyperplanes tangent to
$ Y $), and hence
$ D $ is irreducible.
At a smooth point of
$ D $, the corresponding hyperplane section of
$ Y $ has only one ordinary double point.
We assume that the associated vanishing cycle is not zero in
the cohomology of general hyperplane section
$ X $.
This is equivalent to the non surjectivity of
$ H^{n-1}(Y)\to H^{n-1}(X) $.

A Lefschetz pencil of
$ Y $ is a line
$ \bP^{1} $ in
$ \cP $ intersecting the discriminant
$ D $ at smooth points of
$ D $ (corresponding to hyperplane sections having only one
ordinary double point).
We have a projective morphism
$ \pi : \tY \to \bP^{1} $ such that
$ \tY_{t} := \pi^{-1}(t) $ is the hyperplane section corresponding
$ t \in \bP^{1} \subset \cP $ and
$ \tY $ is the blow-up of
$ Y $ along a smooth closed subvariety
$ Z $ of codimension
$ 2 $ which is the intersection of
$ \tY_{t} $ for any (or two of)
$ t \in \bP^{1} $.

A Lefschetz pencil of hypersurface sections of degree
$ d $ is defined by replacing the embedding of
$ Y $ using
$ \cO_{Y}(d) $ so that a hyperplane section corresponds to a
hypersurface section of degree
$ d $.
Here
$ \cO_{Y}(d) $ for an integer
$ d $ denote the invertible sheaf induced by that on
$ \cP $ as usual.

\medskip\noindent
{\bf 3.2.~Hypersurfaces containing a subvariety.}
Let
$ Y,\cP $ be as in (3.2).
Let
$ E $ be a closed subvariety (which is not necessarily irreducible
nor reduced).
Let
$$
E_{\{i\}} = \{x\in E:\dim T_{x}E=i\}.
$$
Let
$ \cI_{E} $ be the ideal sheaf of
$ E $ in
$ Y $.
Let
$ \delta $ be a positive integer such that
$ \cI_{E}(\delta) $ is generated by global sections.
By [18], [20] (or [21]) we have the following

\medskip\noindent
(3.2.1)\,\,\,
If
$ \dim Y > \max\{\dim E_{\{i\}} + i \} $ and
$ d \ge \delta $, then there is a smooth hypersurface section of
degree
$ d $ containing
$ E $.

\medskip
We have furthermore

\medskip\noindent
(3.2.2)\,\,\,
If
$ \dim Y > \max\{\dim E_{\{i\}} + i \} + 1 $ and
$ d \ge \delta + 1 $, then there is a Lefschetz pencil of hypersurface
sections of degree
$ d $ containing
$ E $.

\medskip
Indeed, we have a pencil such that
$ \tY_{t} $ has at most isolated singularities, because
$ \tY_{t} $ is smooth near the center
$ Z $ which is the intersection of generic two hypersurfaces sections
containing
$ E $, and hence is smooth, see [18], [20] (or [21]).
Note that a local equation of
$ \tY_{t} $ near
$ Z $ is given by
$ f-tg $ if
$ t $ is identified with an appropriate affine coordinate of
$ \bP^{1} $ where
$ f, g $ are global sections of
$ \cI_{E}(d) $ corresponding to smooth hypersurface sections.

To get only ordinary double points, note first that the parameter
space of the hypersurfaces containing
$ E $ is a linear subspace of
$ \cP^{\vee} $.
So it is enough to show that this linear subspace contains a
point of the discriminant
$ D $ corresponding to an ordinary double point.
Thus we have to show that an isolated singularity can be deformed
to ordinary double points by replacing the corresponding section
$ h \in \Gamma(Y,\cI_{E}(d)) $ with
$ h+\sum_{i}t_{i}g_{i} $ where
$ g_{i}\in \Gamma(Y,\cI_{E}(d)) $ and the
$ t_{i}\in \bC $ are general with sufficiently small absolute values.
Since
$ d \ge \delta + 1 $, we see that
$ \Gamma(Y,\cI_{E}(d)) $ generates the
$ 1 $-jets at each point of the complement of
$ E $.
So the assertion follows from the fact that
for a function with an isolated singularity
$ f $, the singularities of
$ \{f + \sum_{i}t_{i}x_{i} = 0\} $ are ordinary double points
if
$ t_{1},\dots,t_{n} $ are general, where
$ x_{1},\dots,x_{n} $ are local coordinates.
(Note that
$ f $ has an ordinary double point if and only if the morphism
defined by
$ (\partial f/\partial x_{1},\dots,\partial f/\partial x_{n}) $ is
locally biholomorphic at this point.)

\medskip\noindent
{\bf 3.3.~Construction.}
For
$ Y,\cP $ be as in (3.1), let
$ i_{Y,\cP} : Y \to \cP $ denote the inclusion.
Assume
$$
i_{Y,\cP}^{*} : H^{j}(\cP)\to H^{j}(Y)\,\,\text{is surjective for any}
\,\, j \ne \dim Y,
\leqno(3.3.1)
$$
where cohomology has coefficients in any field
of characteristic zero.
This condition is satisfied if
$ Y $ is a complete intersection.

Let
$ E_{1} $,
$ E_{2} $ be
$ m $-dimensional irreducible closed subvarieties of
$ Y $ such that
$$
E_{1} \cap E_{2} = \emptyset,\quad\deg E_{1} = \deg E_{2}.
$$
Here
$ \dim Y = n = 2m + s + 1 $ with
$ m \ge 0 $,
$ s \ge 1 $.
Let
$ E = E_{1} \cup E_{2} $.
With the notation of (3.2), assume
$$
d > \delta,\quad
\dim Y > \max\{\dim E_{\{i\}} + i \} + s,
\leqno(3.3.2)
$$
$$
i_{X^{(j)},Y}^{*}:H^{n-j}(Y) \to H^{n-j}(X^{(j)})
\,\, \text{is not surjective for}\,\, j \le s,
\leqno(3.3.3)
$$
where
$ X^{(j)} $ is a general complete intersection of multi degree
$ (d,\dots,d) $ and of codimension
$ j $ in
$ Y $.
(This is equivalent to the condition that the vanishing cycles for
a hypersurface
$ X^{(j)} $ of
$ X^{(j-1)} $ are nonzero.)

Let
$ X $ be a general hypersurface of degree
$ d $ in
$ Y $ containing
$ E $, see (3.2.1).
Let
$ L $ denote the intersection of
$ X $ with a general linear subspace of codimension
$ m + s $ in the projective space.
Then
$ [E_{a}] \,(a = 1,2) $ and
$ c[L\cap X] $ have the same cohomology class in
$ H^{2m+2s}(X) $ for some
$ c \in \bQ $, because
$ \dim H^{2m+2s}(X) = 1 $ by the weak and hard Lefschetz theorems
together with (3.3.1).
Let
$$
\xi_{a}=[E_{a}]-c[L\cap X]\in\CH^{m+s}(X)_{\bQ}\,\,(a=1,2).
$$
These are homologous to zero.
It may be expected that one of them is non torsion,
generalizing an assertion in [24].
More precisely, let
$ S $ be a smooth affine rational variety defined over a finitely
generated subfield
$ k $ of
$ \bC $ and parametrizing the smooth hypersurfaces of degree
$ d $ containing
$ E $ as above so that there is the universal family
$ \cX \to S $ defined over
$ k $ (see [2], [28]).
Assume
$ X $ corresponds to a geometric generic point of
$ S $ with respect to
$ k $, i.e.
$ X $ is the geometric generic fiber for some embedding
$ k(S) \to \bC $.
Let
$$
\xi_{a,\cX} = [E_{a}\mtim_{k}S]-c[L]_{\cX}\in
\CH^{m+s}(\cX)_{\bQ},
$$
where
$ [L]_{\cX} $ is the pull-back of
$ [L] $ by
$ \cX \to Y $.
Since the local system
$ \{H^{2m+2s-j}(\cX_{s})\} $ on
$ S $ is constant for
$ j < s $ and
$ S $ is smooth affine rational, we see that
$ \delta_{S}^{j}(\xi_{a,\cX}) = 0 $ for
$ j < s $.
Then it may be expected that
$ \delta^{s}_{S}(\xi_{a,\cX}) \ne 0 $ for one of
$ a $, where
$ S $ can be replaced by any non empty open subvariety.
We can show this for
$ s = 1 $ as follows.
(For
$ s > 1 $, it may be necessary to assume further conditions on
$ d $, etc.)

\medskip\noindent
{\bf 3.4.~Case
$ s = 1 $.}
Consider a Lefschetz pencil
$ \pi : \tY \to \bP^{1} $ such that
$ \tY_{t} := \pi^{-1}(t) $ for
$ t \in \bP^{1} $ is a hypersurface of degree
$ d $ in
$ Y $ containing
$ E $.
Here
$ \tY $ is the blow-up of
$ Y $ along a smooth closed subvariety
$ Z $, and
$ Z $ is the intersection of
$ \tY_{t} $ for any
$ t \in \bP^{1} $.
Note that
$ \tY_{t} $ has an ordinary double point for
$ t \in \Lambda \subset \bP^{1} $, where
$ \Lambda $ denotes the discriminant, see (3.2.2).

Since
$ Z $ has codimension
$ 2 $ in
$ Y $, we have the isomorphism
$$
H^{n}(\tY) = H^{n}(Y)\oplus H^{n-2}(Z),
\leqno(3.4.1)
$$
so that the cycle class of
$ [E_{a}\times\bP^{1}]-c[L]_{\tY}\in\CH^{m+1}(\tY)_{\bQ} $ in
$ H^{n}(\tY) $ is identified with the difference of the cycle class
$ cl_{Z}(E_{a}) \in H^{n-2}(Z) $ and the cycle class of
$ L $ in
$ H^{n}(Y) $.
Indeed, the injection
$ H^{n-2}(Z) \to H^{n}(\tY) $ in the above direct sum decomposition
is defined by using the projection
$ Z\times\bP^{1}\to Z $ and the closed embedding
$ Z\times\bP^{1}\to \tY $, and the injection
$ H^{n}(Y) \to H^{n}(\tY) $ is the pull-back by
$ \tY \to Y $, see [17].

By assumption, one of the
$ cl_{Z}(E_{a}) $ is not contained in the non primitive part, i.e.
not a multiple of the cohomology class of
the intersection of general hyperplane sections.
Indeed, if both are contained in the non primitive part, then
$ cl_{Z}(E_{1}) = cl_{Z}(E_{2}) $ and this implies the vanishing of
the self intersection number
$ E_{a}\cdot E_{a} $ in
$ Z $.

We will show that the cycle class of
$ [E_{a}\times\bP^{1}]-c[L]_{\tY} $ does not vanish in the
cohomology of
$ \pi^{-1}(U) $ for any non empty open subvariety of
$ \bP^{1} $, in other words, it does not belong to the image of
$ \mopls_{t\in\Lambda'}H_{\tY_{t}}^{n}(\tY) $ where
$ \Lambda' $ is any finite subset of
$ \bP^{1} $ containing
$ \Lambda $.
(Note that the condition for the Lefschetz pencil is generic, and
for any proper closed subvariety of the parameter space, there is a
Lefschetz pencil whose corresponding line is not contained in this
subvariety.)

Thus the assertion is reduced to that
$ \dim H_{\tY_{t}}^{n}(\tY) $ is independent of
$ t \in \bP^{1} $ because this implies that the image of
$ H_{\tY_{t}}^{n}(\tY) \to H^{n}(\tY) $ is independent of
$ t $.
(Note that the Gysin morphism
$ H^{n-2}(\tY_{t}) \to H^{n}(\tY) $ for a general
$ t $ can be identified with the direct sum of the Gysin morphism
$ H^{n-2}(\tY_{t}) \to H^{n}(Y) $ and the restriction morphism
$ H^{n-2}(\tY_{t}) \to H^{n-2}(Z) $ up to a sign, and the image
of the last morphism is the non primitive part by the weak Lefschetz
theorem.)
By duality, this is equivalent to that
$ R^{n}\pi_{*}\bQ_{\tY} $ is a local system on
$ \bP^{1} $.
Then it follows from the assumption that the vanishing cycles
are nonzero, see (3.3.3).

\end{document}